\documentclass[11pt]{article}

\usepackage{amssymb, amsmath, amsthm}

\usepackage{epsfig}
\usepackage{graphicx}
\usepackage{color}
\usepackage{hyperref}
\usepackage[shortalphabetic]{amsrefs}
\usepackage{mathptmx}

\usepackage{fullpage}
\usepackage{pinlabel}
\usepackage{wrapfig}

\theoremstyle{plain}
\newtheorem*{theorem*}{Theorem}
\newtheorem{theorem}{Theorem}

\theoremstyle{remark}

\newtheorem*{remark*}{Remark}

\newtheorem*{challenge}{Challenge}

\newcommand{\R}{{\mathbb R}}

\newcommand{\T}{\mathbb T}
\newcommand{\K}{\mathbb K}

\newcommand{\defn}[1]{\textbf{#1}}
\newcommand{\tun}{\mathfrak{t}}

\newcommand{\spacing}{ \parskip 6.6pt \parindent 0pt}

\spacing

\begin{document}

\title{Two More Proofs that\\ the Kinoshita Graph is Knotted}
\markright{Notes}
\author{Makoto Ozawa and Scott A. Taylor}

\maketitle

\begin{abstract}
The Kinoshita graph is a particular embedding in the 3-sphere of a graph with three edges,  two vertices, and no loops. It has the remarkable property that although the removal of any edge results in an unknotted loop, the Kinoshita graph is itself knotted. We use two classical theorems from knot theory to give two particularly simple proofs that the Kinoshita graph is knotted.
\end{abstract}

\noindent

\begin{figure}[ht!]
\centering
\includegraphics[scale=0.35]{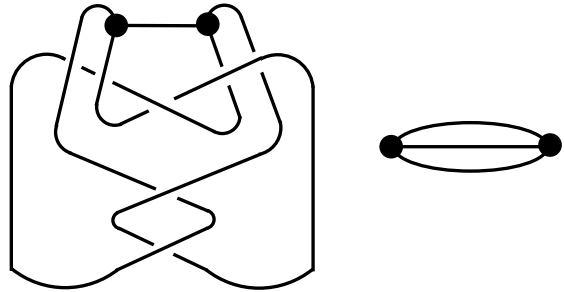}
\caption{On the left is the Kinoshita graph $\K$ and on the right is the trivial theta graph $\T$. For expository ease, we've chosen a diagram of the Kinoshita graph having the property that one edge contains no crossings. Other, prettier, diagrams for the Kinoshita graph can be found by searching online. This diagram is suggested by \cite{Ozawa, SimonWolcott, Wu}.}
\label{Kinoshita}
\end{figure}

The Kinoshita graph $\K$ (left side of Figure \ref{Kinoshita}) is a particular spatial $\theta$-graph; that is, an embedding of a graph with two vertices, three edges, and no loops in the 3-sphere $S^3 = \R^3 \cup \{\infty\}$.  The graph $\K$ has the unusual property that removing any edge from $\K$ results in an unknotted cycle. Of course, the \defn{trivial} spatial $\theta$-graph $\T$, shown on the right of Figure \ref{Kinoshita}, also has this property. For the Kinoshita graph to be interesting, we need to know that there is no continous deformation (techically, an \defn{ambient isotopy}) of $S^3$ taking the Kinoshita graph  $\K$ to the trivial $\theta$-graph $\T$. That is, we need to know the following.
\begin{theorem}\label{thm:nontrivial}
The Kinoshita graph $\K$ is nontrivial.
\end{theorem}

Many proofs of Theorem \ref{thm:nontrivial} are known, including the original one  by Kinoshita \cite{Kinoshita}. See \cite{Litherland, Livingston, MSW, Ozawa12, Scharlemann,  SimonWolcott, Wolcott} for others. Additionally, it is known that $\K$ is prime and hyperbolic. See \cite{CMB, Heard, Ozawa}, \cite[Chapter 3]{Thurston}, \cite[Example 3.3.12]{Thurston2}, and \cite{Wu}. The proofs of these results use a variety of tools from algebraic and geometric topology, including Alexander ideals, branched covers, hyperbolic structures, and cut-and-paste 3-manifold topology. In this note, we present two particularly simple proofs of Theorem \ref{thm:nontrivial}, relying only on classical facts concerning composite knots.

\begin{figure}[ht]
\centering
\includegraphics[scale=.25]{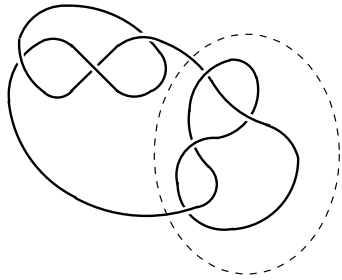}
\caption{A connected sum $\kappa$ of a trefoil knot and a figure 8 knot. The dashed ellipse represents a sphere (the boundary of the 3-balls in the construction described above) in $S^3$ separating the two summands.}
\label{ConnectedSum}
\end{figure}

Recall that $S^3$ is the result of gluing two 3-balls together using any homeomorphism of their boundary. Conversely, the Schoenflies theorem says that every tame 2-sphere in $S^3$ separates $S^3$ into two 3-balls. Given knots $K_1$ and $K_2$ in distinct copies of $S^3$, we can form their \defn{connected sum} as follows. Begin by choosing points $p_1 \in K_1$ and $p_2 \in K_2$. Next, remove and discard open regular neighborhoods of $p_1$ and $p_2$ in the corresponding copies of $S^3$. We are left with two 3-balls $B_1$ and $B_2$. The ball $B_1$  contains a strand which is $K_1 \cap B_1$. Similarly, $B_2$ contains the strand $K_2 \cap B_2$. Choose a homeomorphism $\phi$ between the boundaries of $B_1$ and $B_2$ taking the endpoints of one strand to the endpoints of the other strand. Finally, construct $S^3$ by gluing $B_1$ to $B_2$ using the homeomorphism $\phi$. The union of the strands is a knot $K_1 \# K_2$. There is some ambiguity arising from the choice of $\phi$ and the points $p_1, p_2$, but up to equivalence in $S^3$, at most two different knots can result. (The two possibilities arise from orientation considerations.) Figure \ref{ConnectedSum} shows a connected sum $\kappa$ of a trefoil and a figure 8 knot. A knot that is equivalent (i.e., ambient isotopic to) the connected sum of two nontrivial knots is \defn{composite}; a nontrivial, noncomposite knot is \defn{prime}.

\begin{figure}[ht]
\centering
\includegraphics[scale=.25]{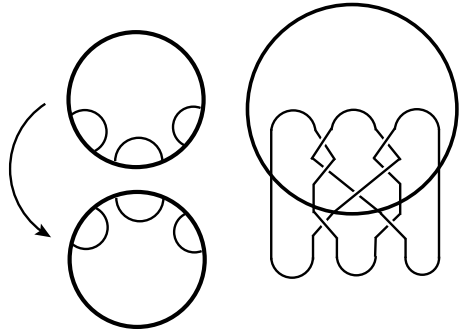}
\caption{On the left is a schematic depiction of gluing two trivial 3-tangles together to produce a knot or link $K$ with $\mathfrak{b}(K) \leq 3$. On the right, we have a knot $K$ in $S^3$ with $\mathfrak{b}(K) \leq 3$. The thick circle denotes a sphere separating the knot into two trivial 3-tangles.}
\label{Bridges}
\end{figure}

Our first proof is probably the simplest possible, though it provides slightly less information about the Kinoshita graph than the second. The knot invariant we'll use for this proof is the bridge number of a knot. Bridge number, like connected sum, is defined using a certain way of gluing two 3-balls together to obtain $S^3$. Consider two 3-balls, each containing the same number $n$ of strands. Unlike in the definition of connected sum, we require that in each 3-ball the strands can be simultaneously isotoped into the boundary of the 3-ball, as in Figure \ref{Bridges}, where $n = 3$. These 3-balls, together with the strands they contain, are called \defn{trivial $n$-tangles}. Trivial $n$-tangles do have diagrams with no crossings, as on the top left and bottom left of Figure \ref{Bridges}; however, they also have diagrams with lots of crossings as on the top right and bottom right of Figure \ref{Bridges}.  We then choose a homeomorphism $\phi$ between the boundaries of the 3-balls, taking the endpoints of the strands in one 3-ball to the endpoints of the strands in the other 3-ball. Gluing the 3-balls together along their boundary using $\phi$ produces the 3-sphere $S^3$ and the union of the strands is a knot or link in $S^3$. Every knot or link in $S^3$ can be obtained this way (for some choice of $n$ and $\phi$). For a given knot or link $K$, the \defn{bridge number} $\mathfrak{b}(K)$ is smallest value of $n$ such that there exists a homeomorphism $\phi$ such that the resulting knot or link is equivalent to $K$. The knot $K$ is the unknot if and only if $\mathfrak{b}(K) = 1$. Schubert \cite{Schubert} proved the following marvelous theorem. (See \cite{Schultens} for a different proof.)

\begin{theorem*}[Schubert]
Suppose that $K_1$ and $K_2$ are knots in $S^3$ and that $K_1 \# K_2$ is any connected sum of them. Then
\[
\mathfrak{b}(K_1 \# K_2) = \mathfrak{b}(K_1) + \mathfrak{b}(K_2) - 1.
\]
\end{theorem*}

One consequence of Schubert's theorem is that if $K$ is a knot with $\mathfrak{b}(K) = 2$, then $K$ is prime.

Suppose now, for a contradiction, that the Kinoshita graph $\mathbb{K}$ is trivial. Then there is an ambient isotopy of $\mathbb{K}$ to the trivial $\theta$-graph $\mathbb{T}$. If $e$ is an edge of $\mathbb{K}$, then this isotopy takes $e$ to an edge $e'$ of $\mathbb{T}$. It also takes an open regular neighborhood of $e$ to an open regular neighborhood of $e'$. The complement of the open regular neighborhood of $e$ is a closed 3-ball $B$. The isotopy takes $B$ to the complement $B'$ of the open regular neighborhood of $e'$. The graph $\mathbb{K}$ intersects $B$ in two strands. Similarly, the graph $\mathbb{T}$ intersects $B'$ in two strands. The isotopy induces a homeomorphism of pairs taking $(B, \mathbb{K} \cap B)$ to $(B', \mathbb{T} \cap B')$.  For each edge of $\mathbb{T}$, it is easy to verify that the complement of an open regular neighborhood of that edge is a trivial $2$-tangle. Thus, $(B, \mathbb{K} \cap B)$ is a trivial 2-tangle. We will use bridge number to show that this is impossible. 

\begin{figure}[ht]
\centering
\includegraphics[scale=.25]{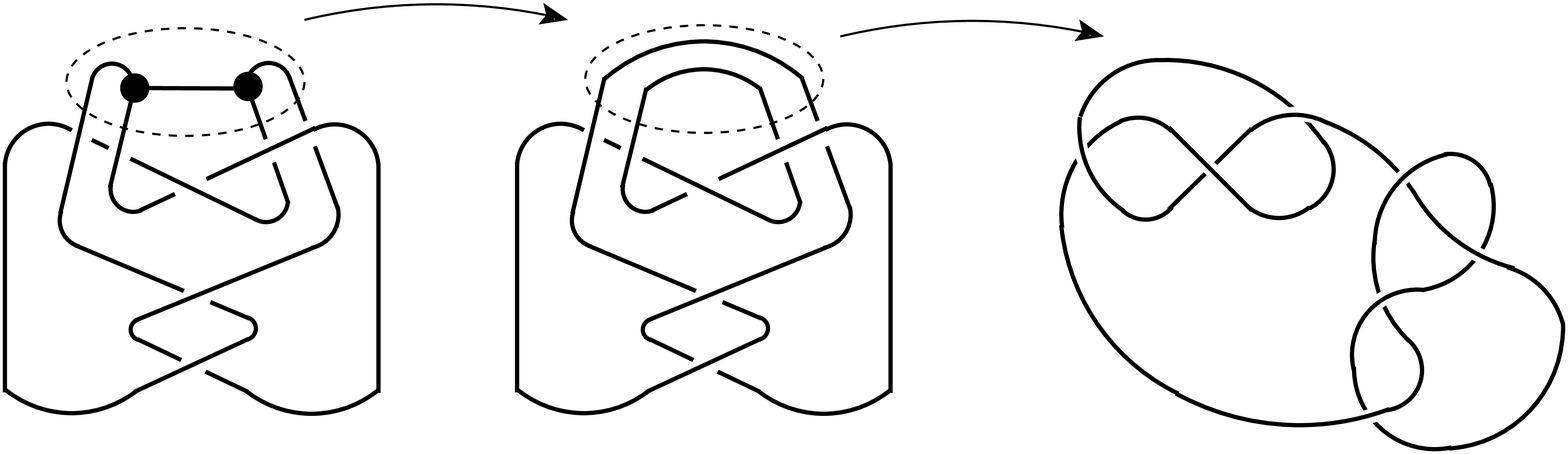}
\caption{On the left the dashed ellipse encloses a 3-ball $N(e)$ which is a regular neighborhood of the edge $e$ of $\K$. The first arrow denotes the action of replacing $(N(e),\K \cap N(e))$ with a certain trivial tangle. The second arrow depicts an isotopy taking the resulting knot to the knot $\kappa$ that is clearly composite.}
\label{KinoshitaTangle}
\end{figure}

Inside the ball $N(e)$ complementary to $B$ replace $\mathbb{K} \cap N(e)$ with a trivial 2-tangle, as in the first step of  Figure \ref{KinoshitaTangle}. We arrive at a knot $K$, which must be a two-bridge knot since it was created by gluing together two trivial 2-tangles. (One tangle is obviously trivial, the other we have shown must be trivial if $\mathbb{K}$ is trivial.) However, as shown in the second step of Figure \ref{KinoshitaTangle}, the knot $K$ is the composite knot $\kappa$! This contradicts Schubert's theorem, and so $\mathbb{K}$ is nontrivial. 

\begin{remark*}
A somewhat more involved argument lets us bypass the use of Schubert's theorem, at the expense of more work. By keeping track of which trivial tangle we place in the ball $N(e)$ and how it is affected by the hypothetical isotopy from $\mathbb{K}$ to $\mathbb{T}$, it is possible to show that the knot $K$ we constructed must be both a torus knot and the connected sum of a trefoil  and a figure 8. This is impossible, as no torus knot is composite. (See \cite{Cromwell} for the idea of how to prove this.)
\end{remark*}

\begin{figure}[ht!]
\centering
\includegraphics[scale=0.35]{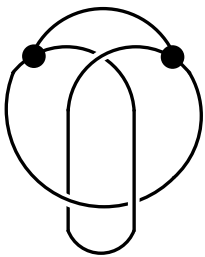}
\caption{An example of a nontrivial spatial $\theta$-graph having handlebody exterior. To see that the graph is nontrivial, note that it contains a knotted cycle (a trefoil knot). Figure \ref{HandlebodyExterior} below shows that the $\theta$-graph  has handlebody exterior, and hence that the trefoil knot has tunnel number one.}
\label{HeegaardSpine}
\end{figure}

We now embark on our second proof, during which we will learn an important fact about the exterior of the Kinoshita graph: it's not a handlebody. A \defn{handlebody} is the result of attaching the ends of solid tubes to 3-balls so as to arrive at a connected, orientable 3-manifold with boundary. More precisely, it is a closed regular neighborhood $N(G)$ of a finite spatial graph $G \subset S^3$, called a \defn{spine} of the handlebody. The \defn{genus} of the handlebody is, by definition, the genus of the bounding surface. Thus, every spatial $\theta$-graph is a spine for a genus 2 handlebody. The trivial $\theta$-graph $\T$ has the property that its exterior $X(\T)$ (the closure of $S^3 \setminus N(G)$) is also a handlebody. Indeed, every spatial graph contained in a tame 2-sphere has this property. The handlebody $X(\T)$ can be constructed using two 3-balls and attaching solid tubes that run through the disk faces. The two handlebodies $N(\T)$ and $X(\T)$ form what is called a \defn{Heegaard splitting} of $S^3$. But there are many other spatial graphs with this property, for instance, the one appearing in Figure \ref{HeegaardSpine}.

\begin{figure}[ht]
\centering
\includegraphics[scale=0.32]{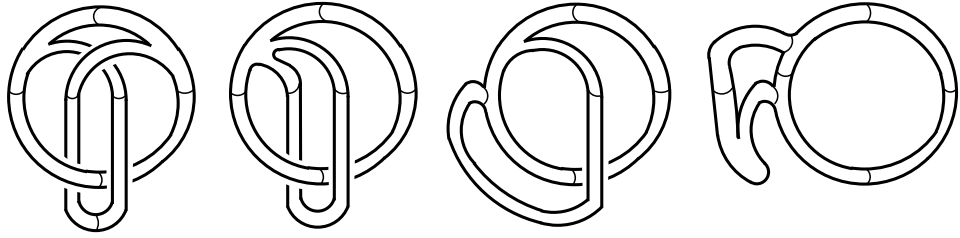}
\caption{Some frames from a movie showing how to perform an ambient isotopy of a regular neighborhood of the graph from Figure \ref{HeegaardSpine} to a regular neighborhood of the trivial $\theta$-graph $\T$.}
\label{HandlebodyExterior}
\end{figure}

An ambient isotopy of a spatial graph $G$ to a spatial graph $G'$ can be extended to an ambient isotopy of $N(G)$ to $N(G')$; however, an ambient isotopy from one handlebody to another need not restrict to an ambient isotopy between two given spines. As no ambient isotopy of a handlebody changes the homeomorphism type of its exterior, this allows us to construct many potentially distinct spatial graphs having homeomorphic exteriors. In particular, every spatial graph $G$ such that $N(G)$ is ambient isotopic to $N(\T)$ has handlebody exterior. Figure \ref{HandlebodyExterior} shows an ambient isotopy of a regular neighborhood of the graph from Figure \ref{HeegaardSpine} to $N(\T)$.  On the other hand, we will show that $\K$ does not have a handlebody exterior, in which case there can be no ambient isotopy taking $N(\K)$ to $N(\T)$, much less one taking $\K$ to $\T$.  This also will show that $\K$ is not equivalent to any graph having handlebody exterior. 

\begin{remark*}
Kinoshita's original proof \cite{Kinoshita} that $\K$ is knotted also proceeds by showing that $\K$ does not have handlebody exterior, but uses new algebraic techniques rather than classical knot invariants.
\end{remark*}

A nontrivial knot $K$ (such as the trefoil knot in Figure \ref{HeegaardSpine})  has \defn{tunnel number one} if it has the property that we can attach an arc in such a way as to create a spatial $\theta$-graph with handlebody exterior. In general, a knot $K$ has tunnel number $\tun(K) = n$ if  $n$ is the smallest integer such that we may attach $n$ arcs to $K$ to arrive at a spatial graph with handlebody exterior \cite{Clark}. Unlike bridge number, tunnel number need not be additive under connected sum (see \cite{Kobayashi, Morimoto95}, for example); however, Norwood did prove the following inequality \cite{Norwood}. See \cite{Sch-Poenaru, TT} for other proofs and \cite{GR, Morimoto, SS} for generalizations.

\begin{theorem*}[Norwood]
If $K_1$ and $K_2$ are nontrivial knots in $S^3$, then $\tun(K_1 \# K_2) \geq 2$.
\end{theorem*}

\begin{figure}[ht]
\centering
\includegraphics[scale=.27]{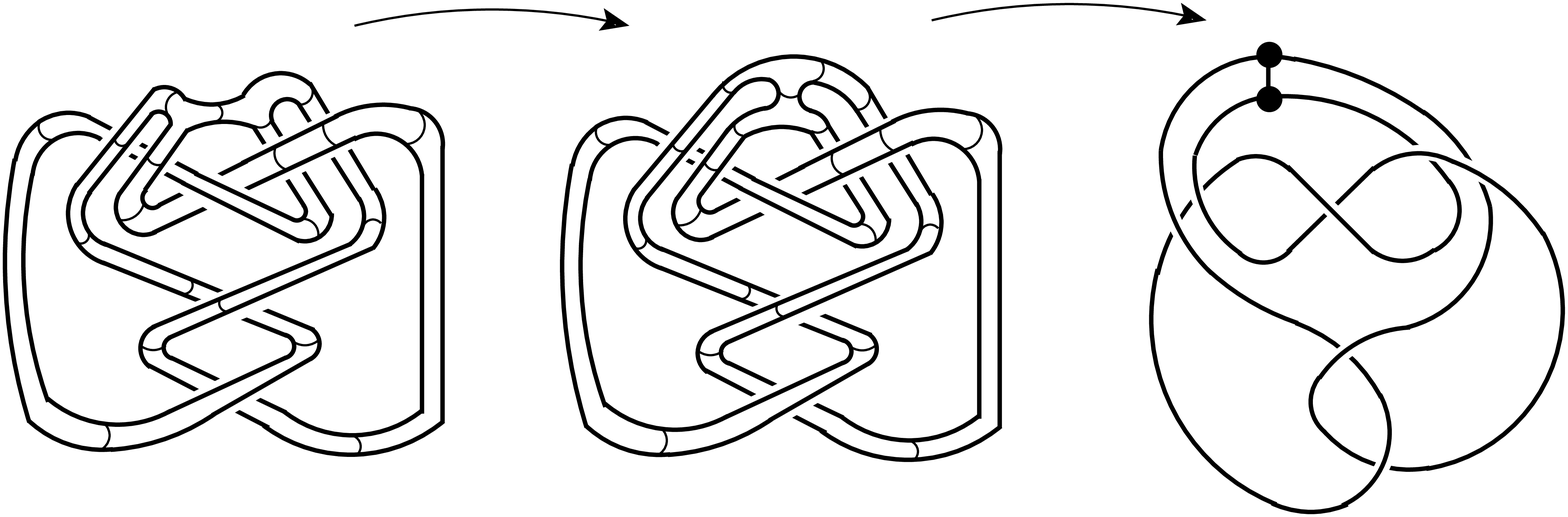}
\caption{The first arrow shows an isotopy of $N(\K)$. The second arrow focuses attention on a spine (different from $\K$) for the handlebody.}
\label{KinoshitaNbhd}
\end{figure}

As a consequence of Norwood's theorem, we see that knots of tunnel number one are prime. Equivalently, the connected sum of two nontrivial knots will never be a cycle in a spatial $\theta$-graph with handlebody exterior. In particular, the composite knot $\kappa$ appearing in Figure \ref{ConnectedSum} and in the middle and right of Figure \ref{KinoshitaTangle} cannot be a cycle in a spine of  spatial $\theta$-graph having handlebody exterior. However, Figure \ref{KinoshitaNbhd} shows an ambient isotopy of a regular neighborhood $N(\K)$ of the Kinoshita graph to a handlebody that is a regular neighborhood of a certain spatial $\theta$-graph $G$. The right side of Figure \ref{KinoshitaNbhd} shows $G$. We see that $G$ contains the knot $\kappa$ as a cycle. Since $\kappa$ is composite, $N(G)$, and hence $N(\K)$, does not have handlebody exterior. This concludes our second proof of Theorem \ref{thm:nontrivial}. We end with a challenge.

\begin{challenge}
Modify the preceding proofs to show the nontriviality of the Kinoshita--Wolcott $\theta_n$-graphs \cite{Wolcott}. (These graphs generalize Suzuki's \cite{Suzuki} generalization of the Kinoshita graph.)
\end{challenge}

%\begin{acknowledgment}{Acknowledgment.}
%We thank Ken Baker for helpful comments.
%\end{acknowledgment}

\vfill\eject

\end{document}